\documentclass[11pt]{amsart}%
\usepackage{amsfonts}
\usepackage{amsmath}
\usepackage{amssymb}
\usepackage{graphicx}%
\providecommand{\U}[1]{\protect\rule{.1in}{.1in}}
\newtheorem{theorem}{Theorem}[section]
\theoremstyle{plain}

\numberwithin{equation}{section}
\begin{document}
\title[Hessian estimate]{Hessian estimates for convex solutions to quadratic Hessian equation}
\author{Matt McGonagle}
\address{University of Washington\\
Department of Mathematics, Box 354350\\
Seattle, WA 98195}
\email{mmcgonal@math.washington.edu}
\author{Chong SONG}
\address{Xiamen University\\
School of Mathematical Sciences\\
Xiamen, 361005, P.R. China}
\email{songchong@xmu.edu.cn}
\author{Yu YUAN}
\address{University of Washington\\
Department of Mathematics, Box 354350\\
Seattle, WA 98195}
\email{yuan@math.washington.edu}

\begin{abstract}
We derive Hessian estimates for convex solutions to quadratic Hessian equation
by compactness argument.
\end{abstract}
\date{\today}
\maketitle

\section{\bigskip Introduction}

In this note, we prove a priori Hessian estimates for convex solutions to the
Hessian equation
\[
\sigma_{k}\left(  D^{2}u\right)  =\sigma_{k}\left(  \lambda\right)
=\sum_{1\leq i_{1}<\cdots<i_{k}\leq n}\lambda_{i_{1}}\cdots\lambda_{i_{k}}=1
\]
with $k=2.$ Here $\lambda_{i}s$ are the eigenvalues of the Hessian $D^{2}u.$

\begin{theorem}
Let $u$ be a smooth solution to $\sigma_{2}\left(  D^{2}u\right)  =1$ on
$B_{R}\left(  0\right)  \subset\mathbb{R}^{n}$ with $D^{2}u\geq\left[
\delta-\sqrt{2/\left[  n\left(  n-1\right)  \right]  }\right]  ~I$ for any
$\delta>0.$ Then%
\[
\left\vert D^{2}u\left(  0\right)  \right\vert \leq g\left(  \left\Vert
Du\right\Vert _{L^{\infty}\left(  B_{R}\left(  0\right)  \right)
}/R,n\right)  ,
\]
where $g\left(  t,n\right)  $ is a finite and positive function for each
positive $t$ and dimension $n.$
\end{theorem}

By Trudinger's gradient estimates for $\sigma_{k}$ equations [T], we can bound
$D^{2}u$ in terms of the solution $u$ in $B_{2R}\left(  0\right)  $ as%
\[
\left\vert D^{2}u\left(  0\right)  \right\vert \leq g\left(  \left\Vert
u\right\Vert _{L^{\infty}\left(  B_{2R}\left(  0\right)  \right)  }%
/R^{2},n\right)  .
\]

Recall any solution to the Laplace equation $\sigma_{1}\left(  D^{2}u\right)
=\bigtriangleup u=1$ enjoys a priori Hessian estimates; yet there are singular
solutions to the three dimensional Monge-Amp\`{e}re equation $\sigma
_{3}\left(  D^{2}u\right)  =\det D^{2}u=1$ by Pogorelov [P], which
automatically generalize to singular solutions to $\sigma_{k}\left(
D^{2}u\right)  =1$ with $k\geq3$ in higher dimensions $n\geq4.$

A long time ago, Heinze [H] achieved a Hessian bound for solutions to equation
$\sigma_{2}\left(  D^{2}u\right)  =1$ in dimension two by two dimension
techniques. Not so long time ago, Hessian bound for $\sigma_{2}\left(
D^{2}u\right)  =1$ in dimension three was obtained via the minimal surface
feature of the \textquotedblleft gradient\textquotedblright\ graph $\left(
x,Du\left(  x\right)  \right)  $ in the joint work with Warren [WY]. Along
this \textquotedblleft integral\textquotedblright\ way, Qiu [Q] has proved
Hessian estimates for solutions to the three dimensional quadratic Hessian
equation with $C^{1,1}$ variable right hand side. Hessian estimates for convex
solutions to general quadratic Hessian equations have also been obtained via a
new pointwise approach by Guan and Qiu [GQ]. Hessian estimates for solutions
to Monge-Amp\`{e}re equation $\sigma_{n}\left(  D^{2}u\right)  =\det D^{2}u=1$
and Hessian equations $\sigma_{k}\left(  D^{2}u\right)  =1$ ($k\geq2$) in
terms of the reciprocal of the difference between solutions and their boundary
values, were derived by Pogorelov [P] and Chou-Wang [CW], respectively, using
Pogorelov's pointwise technique.  Lastly, we also mention Hessian estimates
for solutions to $\sigma_{k}$ as well as $\sigma_{k}/\sigma_{n}$ equations in
terms of certain integrals of the Hessian by Urbas [U1,U2], Bao-Chen-Guan-Ji [BCGJ].

Our argument towards Hessian bound for a semiconvex solution to $\sigma
_{2}\left(  D^{2}\right)  =1$ is through a compactness one. If the Hessian
blows up at the origin, then the slope of the \textquotedblleft
gradient\textquotedblright\ graph $y=Du\left(  x\right)  $ or $\left(
x,Du\left(  x\right)  \right)  $ already blows up everywhere. But one cannot
see this impossible picture directly (Step 1). After a Legendre-Lewy
transformation of the solution $u\left(  x\right)  $ so that the new solution
$\bar{w}\left(  y\right)  $ has bounded nonnegative Hessian; the new
corresponding equation is uniformly elliptic (for any large negative lower
bound for the original Hessian $D^{2}u$); and the new equation is concave
(only under the particular lower Hessian bound $D^{2}u\geq\left[  \delta
-\sqrt{2/\left[  n\left(  n-1\right)  \right]  }\right]  ~I)$ (Step 2). By the
standard Evans-Krylov-Safonov theory, the smooth \textquotedblleft
gradient\textquotedblright\ graph $\left(  D\bar{w}\left(  y\right)
,y\right)  =\left(  x,Du\left(  x\right)  \right)  $ has a zero slope at the
origin (Step 3). Employing the constant rank theorem of Caffarelli-Guan-Ma
[CGM], the zero slope of the \textquotedblleft gradient\textquotedblright%
\ graph $\left(  D\bar{w}\left(  y\right)  ,y\right)  $ propagates everywhere.
The impossible picture of $\left(  x,Du\left(  x\right)  \right)  $ with
infinite slope everywhere becomes clear (Step 4). In passing, we remark that
in dimension two, the solution is already convex and the new equation is just
the Laplace equation, in turn, our compactness argument is elementary.

Finally, the Hessian estimates for general solutions to quadratic Hessian
equation $\sigma_{2}\left(  D^{2}u\right)  =1$ in higher dimension $n\geq4$
still remain an issue to us.

\section{Proof}

We prove Theorem 1.1 by a compactness argument. By scaling $v\left(  x\right)
=u\left(  Rx\right)  /R^{2},$ we assume $R=1.$ Denote $K=\sqrt{2/\left[
n\left(  n-1\right)  \right]  }.$

Step 1. Otherwise, there exist a sequence of solutions $u_{k}$ to $\sigma
_{2}\left(  D^{2}u\right)  =1$ such that%
\begin{align*}
\left\Vert Du_{k}\right\Vert _{L^{\infty}\left(  B_{1}\right)  }  &
\leq\left\Vert Du\right\Vert _{L^{\infty}\left(  B_{1}\right)  },\\
\left(  \delta-K\right)  ~I  &  \leq D^{2}u_{k}.
\end{align*}
and (convergence)%
\[
\left.
\begin{array}
[c]{c}%
\left\vert D^{2}u_{k}\left(  0\right)  \right\vert \rightarrow\infty,\\
Du_{k}\rightarrow Du_{\infty}\ \ \text{in \ }L^{1}\left(  B_{1}^{m} \right)  ,
\end{array}
\right.  \ \ \ \text{as }k\rightarrow\infty,
\]
where $u_{\infty}\in W^{1,1}\left(  B_{1}\right)  $ and $B_{1}^{m}$ denotes
the $m$ dimensional ball $B_{1}^{m}\left(  0\right)  \subset B_{1}=B_{1}%
^{n}\left(  0\right)  \subset\mathbb{R}^{n}$ for all $m=1,\cdots,n.$ The
$L^{1}$ convergence (possibly passing to a sub-convergent sequence, still
denoted by $u_{k};$ we only need $m=1$) comes from the compact Sobolev
embedding for semiconvex $u_{k}\in W^{2,1}\left(  B_{1}^{m}\right)
\hookrightarrow W^{1,1 }\left(  B_{1}^{m}\right)  ,$ as%
\[
\int_{B_{1}^{m}}\left\vert D^{2}u_{k}+K\right\vert dx\leq\int_{B_{1}^{m}%
}\left(  \bigtriangleup u_{k}+nK\right)  dx\leq C\left(  n\right)  \left[
\left\Vert Du\right\Vert _{L^{\infty}\left(  B_{1}\right)  }+1\right]  .
\]

\noindent\textbf{Remark.} Another way to see the above $L^{1}$ convergence for
$H^{n-m}$ almost all $\left(  x_{m+1},\cdots,x_{n}\right)  $ in $B_{1}%
\cap\mathbb{R}^{n-m}$ is the following. From our equation, $\sqrt{\left\vert
\lambda\right\vert ^{2}+2}=\bigtriangleup u_{k},$ then $\int_{B_{1}}\left\vert
D^{2}u_{k}\right\vert dx<\int_{B_{1}}\bigtriangleup udx\leq C\left(  n\right)
\left\Vert Du_{k}\right\Vert _{L^{\infty}\left(  \partial B_{1}\right)  }\leq
C\left(  n\right)  \left\Vert Du\right\Vert _{L^{\infty}\left(  B_{1}\right)
}.$ (We just assume $\left\Vert Du_{k}\right\Vert _{L^{\infty}\left(  \partial
B_{1}\right)  }\leq\left\Vert Du\right\Vert _{L^{\infty}\left(  B_{1}\right)
}$). Now the compact Sobolev embedding coupled with Fubini theorem implies the
\textquotedblleft almost everywhere\textquotedblright\ $L^{1}$ convergence.

\bigskip

Step 2. As in [CY], we make Legendre-Lewy transformation of solutions
$u_{k}\left(  x\right)  $ to solutions $\bar{w}_{k}\left(  y\right)  $ of a
new uniformly elliptic and concave equation with bounded Hessian from both
sides, so that we can extract smoother convergent limit. The Legendre-Lewy
transformation is the Legendre transformation of $w_{k}\left(  x\right)
=u_{k}\left(  x\right)  +K\left\vert x\right\vert ^{2}/2;$ see [L].
Geometrically we re-present the \textquotedblleft gradient\textquotedblright%
\ graph $G:y=Dw_{k}\left(  x\right)  ,$ or $\left(  x,Dw_{k}\left(  x\right)
\right)  \subset\mathbb{R}^{n}\times\mathbb{R}^{n}$ over y-space as another
\textquotedblleft gradient\textquotedblright\ graph in $\mathbb{R}^{2n}.$ Note
that the (canonical) angles between the tangent planes of $G$ and x-space are%
\[
\arctan\left(  \lambda_{i}+K\right)  \in\lbrack\arctan\delta,\frac{\pi}{2})
\]
by the\ semiconvexity assumption $\lambda_{i}\geq\delta-K.$ From this angle
condition and the symmetry of $\left(  D^{2}w\right)  ^{-1},$ it follows that
$G$ can still be represented as a \textquotedblleft gradient\textquotedblright%
\ graph $x=D\bar{w}_{k}\left(  y\right)  ,$ or $\left(  D\bar{w}_{k}\left(
y\right)  ,y\right)  $ over ball $B_{\delta}\left(  0\right)  $ in y-space,
here we may and assume $Du_{k}\left(  0\right)  =0;$ further the (canonical)
angles between the tangent planes of $G$ and y-space are
\[
\arctan\bar{\lambda}_{i}=\frac{\pi}{2}-\arctan\left(  \lambda_{i}+K\right)
\in(0,\frac{\pi}{2}-\arctan\delta],
\]
where $\bar{\lambda}_{i}s$ are the eigenvalues of the Hessian $D^{2}\bar
{w}_{k}.$

Therefore, the function $\bar{w}_{k}\left(  y\right)  $ satisfies in
$B_{\delta}\left(  0\right)  $%
\[
0\leq D^{2}\bar{w}_{k}=\left(  D^{2}u+K\right)  ^{-1}<\frac{1}{\delta}I
\]
and%
\[
q\left(  \bar{\lambda}\left(  D^{2}\bar{w}_{k}\right)  \right)  =\frac
{\sigma_{n-1}\left(  \bar{\lambda}\right)  }{\sigma_{n-2}\left(  \bar{\lambda
}\right)  }=\frac{1}{\left(  n-1\right)  K},
\]
where $\bar{\lambda}_{i}s$ are the eigenvalues of the Hessian $D^{2}\bar
{w}_{k}.$

As proved in [CY, p.661--663], we have

i) the level set $\Gamma=\left\{  \left.  \bar{\lambda}\right\vert \ q\left(
\bar{\lambda}\right)  =1/\left[  \left(  n-1\right)  K\right]  \right\}  $ is convex;

ii) the normal vector $Dq$ of the level set $\Gamma$ is uniformly inside the
positive cone for $\bar{\lambda}_{i}\in\left[  0,\delta^{-1}\right]  ;$

iii) all but one among $\lambda_{i}s$ are uniformly bounded, equivalently all
but one among $\bar{\lambda}_{i}s$ have a uniform positive lower bound, then
$\sigma_{n-2}\left(  \bar{\lambda}\right)  $ has a uniform positive lower bound.

Thus $\bar{w}_{k}$ satisfies a uniformly elliptic and concave equation.

\bigskip

Step 3. By the Evans-Krylov-Safonov theory, there are a subsequence of
$\bar{w}_{k},$ still denoted by $\bar{w}_{k},$ and $\bar{w}_{\infty}\in
C^{2,\alpha}\left(  B_{\delta/2}\left(  0\right)  \right)  $ with
$\alpha=\alpha\left(  n,\delta\right)  >0$ such that%
\[
\bar{w}_{k}\rightarrow\bar{w}_{\infty}\ \ \text{in }C^{2,\alpha}\left(
B_{\delta/2}\left(  0\right)  \right)  ,
\]
then%
\begin{gather*}
q\left(  \bar{\lambda}\left(  D^{2}\bar{w}_{\infty}\right)  \right)
=1/\left[  \left(  n-1\right)  K\right]  ,\\
D^{2}\bar{w}_{\infty}\left(  y\right)  \geq0,\text{ also one and only one
eigenvalue of }D^{2}\bar{w}_{\infty},\ \text{say, }D_{11}\bar{w}_{\infty
}\ \text{is }0\ \text{at }0.
\end{gather*}

\bigskip

Step 4. By the constant rank theorem of Caffarelli-Guan-Ma [CGM, Theorem 1.1
and Remark 1.7] (which leads to a qualitative lower Hessian bound for concave
equations), $D_{11}\bar{w}_{\infty}\left(  y\right)  \equiv0$ in a
neighborhood of $0.$ Restrict to $\left(  x_{1},y_{1}\right)  $ space, the
\textquotedblleft gradient\textquotedblright\ graph of $\left(  D\bar
{w}_{\infty}\left(  y\right)  ,y\right)  $ takes the form%
\[
\left(  D_{1}\bar{w}_{\infty}\left(  y_{1},y^{\prime}\right)  ,y_{1}\right)
=\left(  c,y_{1}\right)  =\left(  x_{1},D_{1}u_{\infty}\left(  x_{1}%
,x^{\prime}\right)  +Kx_{1} \right)  \ \ \ \text{near }\left(  0,0\right)  .
\]
This is impossible, as $\left(  x_{1},D_{1}u_{\infty}\left(  x_{1},x^{\prime
}\right)  +Kx_{1} \right)  $ is an $L^{1}$ graph (for almost all $x^{\prime
}\in\mathbb{R}^{n-1}$ without using the semiconvexity assumption) from Step 1.

\bigskip

\textbf{Acknowledgments. } Part of this work was carried out when the second
author was visiting University of Washington, which was supported by the China
Scholarship Council. The second author is partially supported by the
Fundamental Research Funds for the Central University (No. 20720170009). The
third author is partially supported by an NSF grant.


\bigskip

\begin{thebibliography}{9999}                                                                                             %


\bibitem[BCGJ]{BCGJ}Bao, Jiguang; Chen, Jingyi; Guan, Bo; Ji, Min
\emph{Liouville property and regularity of a Hessian quotient equation.} Amer.
J. Math. \textbf{125} (2003), no. 2, 301--316.

\bibitem[CGM]{CGM}Caffarelli, Luis; Guan, Pengfei; Ma, Xi-Nan \emph{A constant
rank theorem for solutions of fully nonlinear elliptic equations.} Comm. Pure
Appl. Math. \textbf{60} (2007), no. 12, 1769--1791.

\bibitem[CY]{CY}Chang, Sun-Yung Alice; Yuan, Yu \emph{A Liouville problem for
the sigma-2 equation.} Discrete Contin. Dyn. Syst. \textbf{28} (2010), no. 2, 659--664.

\bibitem[CW]{CW}Chou, Kai-Seng; Wang, Xu-Jia, \emph{A variational theory of
the Hessian equation.} Comm. Pure Appl. Math. \textbf{54} (2001), 1029--1064.

\bibitem[GQ]{GQ}Guan, Pengfei; Qiu, Guohuan, \emph{Interior }$C^{2}$\emph{
regularity of convex solutions to prescribing scalar curvature equations},
preprint (2016).

\bibitem[H]{H}Heinz, Erhard, \emph{On elliptic Monge-Amp\`{e}re equations and
Weyl's embedding problem}. J. Analyse Math. \textbf{7} (1959) 1--52.

\bibitem[L]{L}Lewy, Hans, \emph{A priori limitations for solutions of
Monge-Amp\`{e}re equations.} II. Trans. Amer. Math. Soc. \textbf{41} (1937), 365--374.

\bibitem[P]{P}Pogorelov, Aleksei Vasil'evich, \emph{The Minkowski
Multidimensional Problem.} Translated from the Russian by Vladimir Oliker.
Introduction by Louis Nirenberg. Scripta Series in Mathematics. V. H. Winston
\& Sons, Washington, D.C.; Halsted Press [John Wiley \& Sons], New
York-Toronto-London, 1978.

\bibitem[Q]{Q}Qiu, Guohuan, Interior Hessian estimates for sigma-2 equations,
preprint (2017).

\bibitem[T]{T}Trudinger, Neil S., \emph{Weak solutions of Hessian equations.}
Comm. Partial Differential Equations \textbf{22} (1997), no. 7-8, 1251--1261.

\bibitem[U1]{U1}Urbas, John, \emph{Some interior regularity results for
solutions of Hessian equations.} Calc. Var. Partial Differential Equations
\textbf{11} (2000), 1--31.

\bibitem[U2]{U2}Urbas, John, \emph{An interior second derivative bound for
solutions of Hessian equations.} Calc. Var. Partial Differential Equations
\textbf{12} (2001), 417--431.

\bibitem[WY]{WY}Warren, Micah; Yuan, Yu, \emph{Hessian estimates for the
sigma-2 equation in dimension three.} Comm. Pure Appl. Math. \textbf{62}
(2009), no. 3, 305--321.
\end{thebibliography}
\end{document}